\documentclass[12pt,oneside,reqno]{amsart}
\usepackage{graphicx}
\usepackage{mathrsfs}
\usepackage{stmaryrd}
\usepackage{amsfonts}
\usepackage{enumerate,amsmath,amssymb,amsthm}
\newcommand{\dif}{\mathrm{d}}

\newcommand{\be}{\begin{eqnarray}}
\newcommand{\ee}{\end{eqnarray}}
\newcommand{\ce}{\begin{eqnarray*}}
\newcommand{\de}{\end{eqnarray*}}
\newtheorem{theorem}{Theorem}[section]
\newtheorem{lemma}[theorem]{Lemma}
\newtheorem{remark}[theorem]{Remark}
\newtheorem{definition}[theorem]{Definition}
\newtheorem{proposition}[theorem]{Proposition}
\newtheorem{Example}[theorem]{Example}
\newtheorem{corollary}[theorem]{Corollary}

\def\[{{\Big[}}
\def\]{{\Big]}}
\def\<{{\langle}}
\def\>{{\rangle}}
\def\({{\Big(}}
\def\){{\Big)}}

\def\tr{{\rm tr}}

\def\no{\nonumber}
\def\bt{\begin{theorem}}
\def\et{\end{theorem}}
\def\bl{\begin{lemma}}
\def\el{\end{lemma}}
\def\br{\begin{remark}}
\def\er{\end{remark}}
\def\bx{\begin{Example}}
\def\ex{\end{Example}}
\def\bd{\begin{definition}}
\def\ed{\end{definition}}
\def\bp{\begin{proposition}}
\def\ep{\end{proposition}}
\def\bc{\begin{corollary}}
\def\ec{\end{corollary}}

\def\cL{{\mathcal L}}
\def\cM{{\mathcal M}}
\def\cN{{\mathcal N}}

\def\mE{{\mathbb E}}

\def\mP{{\mathbb P}}
\def\mQ{{\mathbb Q}}
\def\mR{{\mathbb R}}

\def\mU{{\mathbb U}}

\def\sB{{\mathscr B}}
\def\sC{{\mathscr C}}

\def\sF{{\mathscr F}}

\def\sU{{\mathscr U}}

\def\geq{\geqslant}
\def\leq{\leqslant}

\begin{document}

\allowdisplaybreaks

\title{Path independence of the additive functionals for McKean-Vlasov stochastic differential equations with jumps*}

\author{Huijie Qiao$^{1,2}$ and Jiang-Lun Wu$^3$}



\thanks{*This work was partly supported by NSF of China (No. 11001051, 11371352, 11671083) and China Scholarship Council under Grant No. 201906095034.}


\date{}

\dedicatory{1. School of Mathematics,
Southeast University\\
Nanjing, Jiangsu 211189,  China\\
2. Department of Mathematics, University of Illinois at
Urbana-Champaign\\
Urbana, IL 61801, USA\\
hjqiaogean@seu.edu.cn\\
3. Department of Mathematics, Computational Foundry, Swansea University\\
Bay Campus, Swansea SA1 8EN, UK\\
j.l.wu@swansea.ac.uk}

\begin{abstract} 

In this article, the path independent property of additive functionals of McKean-Vlasov stochastic 
differential equations with jumps is characterised by nonlinear partial integro-differential equations involving 
$L$-derivatives with respect to probability measures introduced by P.-L. Lions. Our result extends the recent 
work \cite{rw} by Ren and Wang where their concerned McKean-Vlasov stochastic differential 
equations are driven by Brownian motions. 

\end{abstract}

\maketitle \rm

\section{Introduction}

Since the seminal work \cite{M1,M2}, there have been substantial interests to study McKean-Vlasov stochastic 
differential equations, which are stochastic differential equations whose coefficients depend on the law of the solution,  
and are also referred as mean-field stochastic differential equations, see, e.g., \cite{Lion} and most recently \cite{BLPR,L} 
(and references therein). Very recently,  Ren and Wang \cite{rw} explored an interesting result characterising 
the path independent additive functionals of McKean-Vlasov stochastic differential equations driven by Brownian 
motion by space-distribution partial differential equations, which extends the earlier work \cite{twwy,wy} on this direction.  
  
The object of this paper is to extend \cite{rw} to the same type of equations driven by compensated Poisson martingale 
measures (and Brownian motion). We aim to characterise the path-independence of additive functionals of McKean-Vlasov 
stochastic differential equations with jumps by certain partial integro-differential equations involving $L$-derivatives with 
respect to probability measures, following our previous work \cite{qw1,qw2} where therein stochastic differential equations 
with jumps in finite and infinite dimensions were studied, respectively.  Let us also mention further interesting work \cite{ry,lxz}, 
where characterisation theorems for the path independence of additive functionals of stochastic differential equations driven
by $G$-Brownian motion as well as for stochastic differential equations driven by Brownian motion with non-Markovian 
coefficients (i.e. random coefficients) are established, respectively. It is of course very interesting to extend the two cases 
to the situation of the concerned equations with jumps, which we will consider in our forthcoming work. 

It is worthwhile to mention our results. We prove the It\^o formula for McKean-Vlasov stochastic differential equations. And its proof is more simple than that in \cite{HL} and \cite{L}. Moreover, it does not contain any abstract probability space. Therefore, it is more applicable. Besides, we compare our main results with that in \cite{rw} and \cite{qw1}. And we find that it is indeed more general.

The rest of the paper is organized as follows. In the next section, we will set up the framework for introducing the 
McKean-Vlasov stochastic differential equations. In Section 3, we will first derive the It\^o formula for the solutions of 
our concerned McKean-Vlasov stochastic differential equations (and the proof is given in the Appendix at the end 
of our paper), prove our characterisation theorem, analyze some special cases and finally compare our results with some known results. 
 
\section{Preliminary}\label{pre}

\subsection{Notation}\label{nn}

In the subsection, we introduce notation used in the sequel. 

For convenience, we shall use $\mid\cdot\mid$ and $\parallel\cdot\parallel$  for norms of vectors and matrices, respectively. Furthermore, let $\langle\cdot$ , $\cdot\rangle$ denote the scalar product in $\mR^d$. Let $A^*$ denote the transpose of the matrix $A$.

Let $\sB(\mR^d)$ be the Borel $\sigma$-algebra on $\mR^d$ and $\cM({\mR^d})$ be the space of all probability measures defined on $\sB(\mR^d)$ carrying the usual topology of weak convergence. Let $\cM_2(\mR^d)$ be the collection of all the probability measures $\mu$ on $\sB(\mR^d)$ satisfying
\ce
\mu(|\cdot|^2):=\int_{\mR^d}\mid{x}\mid^2\,\mu(\dif x)<\infty.
\de
We put on $\cM_2(\mR^d)$ a topology induced by the
following metric:
\ce
\rho^2(\mu_1,\mu_2):=\inf_{\pi\in\sC(\mu_1, \mu_2)}\int_{\mR^d\times\mR^d}|x-y|^2\pi(\dif x, \dif y), \quad \mu_1, \mu_2\in\cM_2(\mR^d),
 \de
where $\sC(\mu_1, \mu_2)$ denotes the set of  all the probability measures whose marginal distributions are $\mu_1, \mu_2$, respectively. Thus, $(\cM_2(\mR^d), \rho)$ is a Polish space.

\subsection{McKean-Vlasov stochastic differential equations with jumps}\label{smvej}

In the subsection, we introduce McKean-Vlasov stochastic differential equations with jumps and  path-independence for a type of additive functionals.

Let $(\Omega,\sF,\mP)$ be a complete, filtered probability space. Let $(B_t)_{t\geq 0}$ be a $m$-dimensional Brownian motion. Let $(\mU,\|\cdot\|_{\mU})$ be a finite dimensional normed space with its Borel $\sigma$-algebra $\mathscr{U}$.  Let $\nu$ be a $\sigma$-finite measure defined on $(\mU,\mathscr{U})$. We fix $\mU_0=\{\|u\|_{\mU}\leq\alpha\}$, where $\alpha>0$ is a constant, with $\nu(\mU\setminus\mU_0)<\infty$ and $\int_{\mU_0}\|u\|_{\mU}^2\,\nu(\dif u)<\infty$. Let $N(\dif t, \dif u)$ be an integer-valued Poisson random measure  on $(\Omega,\sF,\mP)$ with intensity $\mE N(\dif t, \dif u)=\dif t\nu(\dif u)$. Denote
$$
\tilde{N}(\dif t,\dif u):=N(\dif t,\dif u)-\dif t\nu(\dif u),
$$
that is, $\tilde{N}(\dif t,\dif u)$ stands for the compensated martingale measure of $N(\dif t, \dif u)$. Moreover, $B_{\cdot}$ and $N(\dif t, \dif u)$ are mutually independent. Fix $T>0$ and let $\{\sF_t\}_{t\in[0,T]}$ be the filtration generated by $(B_t)_{t\geq 0}$ and 
$N(\dif t, \dif u)$, and augmented by a $\sigma$-field $\sF^0$, i.e.,
\ce
&&\sF^0_t:=\sigma\{B_s, N((0,s], A): 0\leq s\leq t, A\in\sU\},\\
&&\sF_t:=(\cap_{s>t}\sF^0_s)\vee\sF^0, \quad t\in[0,T],
\de
where $\sF^0\subset\sF$ has the following properties:

(i) $(B_t)_{t\geq 0}$ and $N(\dif t, \dif u)$ are independent of $\sF^0$;

(ii) $\cM_2(\mR^d)=\{\mP\circ \xi^{-1}, \xi\in L^2(\sF^0; \mR^d)\}$;

(iii) $\sF^0\supset \cN$ and $\cN$ is the collection of all $\mP$-null sets.

Now, consider the following McKean-Vlasov stochastic differential equation with jumps on $\mR^d$:
\be
\dif X_t=b(t,X_t,\cL_{X_t})\dif t+\sigma(t,X_t,\cL_{X_t})\dif B_t+\int_{\mU_0}f(t,X_{t-},\cL_{X_t},u)\tilde{N}(\dif t, \dif u), \qquad t\in[0,T],
\label{Eq1}
\ee
where $\cL_{X_t}$ denotes the distribution of $X_t$ under $\mP$. Here the coefficients $b: [0,T]\times\mR^d\times\cM_2(\mR^d)\mapsto\mR^d$, $\sigma: [0,T]\times\mR^d\times\cM_2(\mR^d)\mapsto\mR^{d\times m}$ and
$f: [0,T]\times\mR^d\times\cM_2(\mR^d)\times\mU_0\mapsto\mR^d$ are all Borel measurable. We assume:
\begin{enumerate}[(${\bf H}_{b, \sigma}$)]
\item 
There exists an increasing function $L_1: [0, \infty)\mapsto (0, \infty)$ such that for $t\in[0,T]$, $x_1, x_2\in\mR^d$, $\mu_1, \mu_2\in\cM_2(\mR^d)$,
\ce
|b(t, x_1, \mu_1)-b(t, x_2, \mu_2)|+\|\sigma(t, x_1, \mu_1)-\sigma(t, x_2, \mu_2)\|\leq L_1(t)(|x_1-x_2|+\rho(\mu_1, \mu_2)),
\de
and for $t\in[0,T]$,
\ce
|b(t,0,\delta_0)|+\|\sigma(t,0,\delta_0)\|\leq L_1(t),
\de
where $\delta_0$ is the Dirac measure in $0$.
\end{enumerate}
\begin{enumerate}[(${\bf H}_{f}$)]
\item There exists an increasing function $L_2: [0, \infty)\mapsto (0, \infty)$ such that for $t\in[0,T]$, $x_1, x_2\in\mR^d$, $\mu_1, \mu_2\in\cM_2(\mR^d)$, $u\in\mU_0$,
\ce
|f(t, x_1, \mu_1, u)-f(t, x_2, \mu_2, u)|\leq L_2(t)\|u\|_{\mU}(|x_1-x_2|+\rho(\mu_1, \mu_2)),
\de
and for $t\in[0,T]$, $u\in\mU_0$,
\ce
|f(t,0,\delta_0,u)|\leq L_2(t)\|u\|_{\mU}.
\de
\end{enumerate}

\br\label{lingro}
By (${\bf H}_{b, \sigma}$) and (${\bf H}_{f}$), it holds that for $t\in[0,T]$, $x\in\mR^d$, $\mu\in\cM_2(\mR^d)$
\ce
&&|b(t,x,\mu)|^2+\|\sigma(t,x,\mu)\|^2\leq C\(1+|x|^2+\mu(|\cdot|^2)\),\\
&&|f(t,x,\mu,u)|^2\leq C\|u\|^2_{\mU}\(1+|x|^2+\mu(|\cdot|^2)\).
\de
\er

Under (${\bf H}_{b,\sigma}$) (${\bf H}_f$), following up the proof of \cite[Theorem 3.1, Page 7]{HL}, we can obtain that for any $s\in[0, T]$ and $X_s\in L^2(\Omega, \sF_s,\mP; \mR^d)$, Eq.(\ref{Eq1}) has a unique solution $(X_t)_{t\geq s}$ with 
\be
\mE\(\sup\limits_{t\in[s,T]}|X_t|^2\)<\infty.
\label{xbou}
\ee

And then we introduce the following additive functional
\be
F_{s,t}&:=&\int_s^t g_1(r,X_r, \cL_{X_r})\dif r+\int_s^t \<g_2(r,X_r, \cL_{X_r}), \dif B_r\>+\int_s^t\int_{\mU_0}g_3(r,X_{r-}, \cL_{X_r}, u)\tilde{N}(\dif r, \dif u)\no\\
&&+\int_s^t\int_{\mU_0}g_4(r,X_r, \cL_{X_r}, u)\nu(\dif u)\dif r, \quad 0\leq s< t\leq T,
\label{addfun}
\ee
where 
\ce
&&g_1: [0,T]\times\mR^d\times\cM_2(\mR^d)\mapsto\mR, \qquad g_2: [0,T]\times\mR^d\times\cM_2(\mR^d)\mapsto\mR^m,\\
&&g_3:[0,T]\times\mR^d\times\cM_2(\mR^d)\times\mU_0\mapsto\mR, g_4:[0,T]\times\mR^d\times\cM_2(\mR^d)\times\mU_0\mapsto\mR,
\de
are Borel measurable, $g_1(t,x,\mu), g_2(t,x,\mu), g_3(t,x,\mu,u)$ are continuous in $(t,x,\mu)$ and $\int_{\mU_0}g_4(t,x,\mu,u)\nu(\dif u)$ is continuous in $(t,x,\mu)$, so that $F_{s,t}$ is a well-defined semimartingale.

\bd\label{pathinde}
The additive functional $F_{s,t}$ is called path independent, if there exists a function 
$$
V: [0,T]\times\mR^d\times\cM_2(\mR^d)\mapsto\mR,
$$
such that for any $s\in[0, T]$ and $X_s\in L^2(\Omega, \sF_s,\mP; \mR^d)$, the solution $(X_t)_{t\in[s,T]}$ of Eq.(\ref{Eq1}) satisfies
\be
F_{s,t}=V(t,X_t,\cL_{X_t})-V(s,X_s,\cL_{X_s}).
\label{defi}
\ee
\ed

\subsection{L-derivative for functions on $\cM_2(\mR^d)$}\label{lde} In the subsection we recall the definition of L-derivative for functions on $\cM_2(\mR^d)$. And the definition was first introduced by Lions \cite{Lion}. Moreover, he used some abstract probability spaces to describe the L-derivatives. Here, for the convenience to understand the definition, we apply a straight way to state it (\cite{rw}). Let $I$ be the identity map on $\mR^d$. For $\mu\in\cM_2(\mR^d)$ and $\phi\in L^2(\mR^d, \sB(\mR^d), \mu;\mR^d)$, $\mu(\phi):=\int_{\mR^d}\phi(x)\mu(\dif x)$. Moreover, by simple calculation, it holds that $\mu\circ(I+\phi)^{-1}\in\cM_2(\mR^d)$.

\bd\label{lderi}
(i) A function $h: \cM_2(\mR^d)\mapsto\mR$ is called L-differentiable at $\mu\in\cM_2(\mR^d)$, if the functional 
$$
L^2(\mR^d, \sB(\mR^d), \mu;\mR^d)\ni\phi\mapsto h(\mu\circ(I+\phi)^{-1})
$$
is Fr\'echet differentiable at $0\in L^2(\mR^d, \sB(\mR^d), \mu;\mR^d)$; that is, there exists a unique $\xi\in L^2(\mR^d, \sB(\mR^d), \mu;\mR^d)$ such that 
$$
\lim\limits_{\mu(|\phi|^2)\rightarrow 0}\frac{h(\mu\circ(I+\phi)^{-1})-h(\mu)-\mu(\<\xi,\phi\>)}{\sqrt{\mu(|\phi|^2)}}=0.
$$
In the case, we denote $\partial_{\mu}h(\mu)=\xi$ and call it the L-derivative of $h$ at $\mu$.

(ii) A function $h: \cM_2(\mR^d)\mapsto\mR$ is called L-differentiable on $\cM_2(\mR^d)$ if  L-derivative $\partial_{\mu}h(\mu)$ exists for all $\mu\in\cM_2(\mR^d)$.

(iii) By the same way, $\partial^2_\mu h(\mu)(y,y')$ for $y, y'\in\mR^d$ can be defined.
\ed

Next, we introduce some related spaces.

\bd\label{space1}
 The function $h$ is said to be in $C^2(\cM_2(\mR^d))$, if $\partial_\mu h$ is continuous, for any $\mu\in \cM_2(\mR^d)$, $\partial_\mu h(\mu)(\cdot)$ is differentiable, and its derivative $\partial_y\partial_\mu h:\cM_2(\mR^d)\times\mR^d\rightarrow\mR^d\otimes\mR^d$ is continuous, and for any $y\in\mR^d$, $\partial_\mu h(\cdot)(y)$ is differentiable, and its derivative $\partial^2_\mu h:\cM_2(\mR^d)\times\mR^d\times\mR^d\rightarrow\mR^d\otimes\mR^d$ is continuous.
\ed

\bd\label{space2}
(i) The function $h: [0,T]\times\mR^d\times\cM_2(\mR^d)\mapsto\mR$ is said to be in $C^{1,2,2}([0,T]\times\mR^d\times\cM_2(\mR^d))$, if $h(t,x,\mu)$ is $C^1$ in $t\in[0,T]$, $C^2$ in $x\in\mR^d$ and $\mu\in\cM_2(\mR^d)$ respectively, and its derivatives 
$$
\partial_t h(t,x,\mu), \partial_x h(t,x,\mu), \partial^2_x h(t,x,\mu), \partial_\mu h(t,x,\mu)(y),  \partial_y\partial_\mu h(t,x,\mu)(y), \partial^2_\mu h(t,x,\mu)(y, y')
$$ 
are jointly continuous in the corresponding variable family $(t,x,\mu)$, $(t,x,\mu,y)$ or $(t,x,\mu,y, y')$. 

(ii) The function $h: [0,T]\times\mR^d\times\cM_2(\mR^d)\mapsto\mR$ is said to be in $C_b^{1,2,2}([0,T]\times\mR^d\times\cM_2(\mR^d))$, if $h\in C^{1,2,2}([0,T]\times\mR^d\times\cM_2(\mR^d))$ and all its derivatives are uniformly bounded on $[0,T]\times\mR^d\times\cM_2(\mR^d)$. If $h\in C^{1,2,2}([0,T]\times\mR^d\times\cM_2(\mR^d))$ or $h\in C_b^{1,2,2}([0,T]\times\mR^d\times\cM_2(\mR^d))$ and $h$ is independent of $t$, we write $h\in C^{2,2}(\mR^d\times\cM_2(\mR^d))$ or $h\in C_b^{2,2}(\mR^d\times\cM_2(\mR^d))$.

(iii) The function $h: [0,T]\times\mR^d\times\cM_2(\mR^d)\mapsto\mR$ is said to be in $C_b^{1,2,2;1}([0,T]\times\mR^d\times\cM_2(\mR^d))$, if $h\in C_b^{1,2,2}([0,T]\times\mR^d\times\cM_2(\mR^d))$ and all its derivatives are Lipschitz continuous. In addition, if $h$ is independent of $t$, we write $h\in C_b^{2,2;1}(\mR^d\times\cM_2(\mR^d))$.
\ed

\section{Main results and related analysis}\label{main}

In the section, we state and prove the main results, analyze some special cases and compare our results with some known results. 

\subsection{Main results and their proofs}
In the subsection, we state and prove the main results.

First of all, we prove the It\^o formula which is an important tool in our following proofs.

\bp\label{itoii}{\bf (The It\^o formula)}
Suppose that (${\bf H}_{b,\sigma}$) (${\bf H}_f$) hold. Then, if $h$ belongs to $C_b^{1,2,2}([0,T]\times\mR^d\times\cM_2(\mR^d))$ and all the derivatives of $h$ in $(t,x,\mu)$ are uniformly continuous, it holds that for $t\geq 0$,
\be
\dif h(t,X_t,\cL_{X_t})&=&(\partial_t+{\bf L}_{b,\sigma,f})h(t,X_t,\cL_{X_t})\dif t+\<(\sigma^*\partial_x h)(t,X_t,\cL_{X_t}), \dif B_t\>\no\\
&&+\int_{\mU_0}\left[h(t, X_{t-}+f(t,X_{t-},\cL_{X_t},u),\cL_{X_t})-h(t, X_{t-},\cL_{X_t})\right]\tilde{N}(\dif t, \dif u), \no\\
\label{itoeq}
\ee
where 
\be
{\bf L}_{b,\sigma, f}h(t, x,\mu)&:=&\<b, \partial_x h\>(t, x,\mu)+\frac{1}{2}\tr\((\sigma\sigma^*)\partial_x^2 h\)(t, x,\mu)\no\\
&&+\int_{\mR^d}\<b(t, y,\mu), (\partial_\mu h)(t, x,\mu)(y)\>\mu(\dif y)\no\\
&&+\frac{1}{2}\int_{\mR^d}\tr\((\sigma\sigma^*)(t, y,\mu)\partial_y\partial_\mu h(t, x,\mu)(y)\)\mu(\dif y) \no\\
&&+\int_{\mU_0}\bigg[h\(t, x+f(t,x,\mu,u),\mu\)-h(t, x,\mu)\no\\
&&\qquad\qquad\qquad -\<f(t,x,\mu,u), \partial_x h(t, x,\mu)\>\bigg]\nu(\dif u)\no\\
&&+\int_{\mU_0}\int_0^1\int_{\mR^d}\bigg<\partial_\mu h(t,x,\mu)\(y+\eta f(t,y,\mu,u)\)\no\\
&&\qquad\qquad\qquad -\partial_\mu h(t,x,\mu)(y), f(t,y,\mu,u)\bigg>\mu(\dif y)\dif \eta \nu(\dif u).
\label{gene}
\ee
\ep

Because the proof of Proposition \ref{itoii} is too long, we place it to the Appendix so as to make the context more compact. Now, it is the position to state and prove our main result.
\bt\label{funthe}
Assume that $b, \sigma, f$ satisfy (${\bf H}_{b,\sigma}$) (${\bf H}_f$). Then if $V$ belongs to $C_b^{1,2,2}([0,T]\times\mR^d\times\cM_2(\mR^d))$ and all the derivatives of $V$ in $(t,x,\mu)$ are uniformly continuous, $g_1\in C([0,T]\times\mR^d\times\cM_2(\mR^d)\mapsto\mR)$, $g_2\in C([0,T]\times\mR^d\times\cM_2(\mR^d)\mapsto\mR^m)$, $g_3(\cdot,\cdot, \cdot, u)\in C([0,T]\times\mR^d\times\cM_2(\mR^d)\mapsto\mR)$ and $\int_{\mU_0}g_4(\cdot,\cdot, \cdot, u)\nu(\dif u)\in C([0,T]\times\mR^d\times\cM_2(\mR^d)\mapsto\mR)$, $F_{s,t}$ is path independent in the sense of (\ref{defi}) if and only if $(V, g_1, g_2, g_3, g_4)$ satisfies the partial integral-differential equation
\be\left\{\begin{array}{ll}
(\partial_t+{\bf L}_{b,\sigma,f})V(t,x,\mu)=g_1(t, x,\mu)+\int_{\mU_0}g_4(t, x,\mu,u)\nu(\dif u),\\
(\sigma^*\partial_x V)(t, x,\mu)=g_2(t, x,\mu),\\
V\(t, x+f(t,x,\mu,u),\mu\)-V(t, x,\mu)=g_3(t, x,\mu,u), \\
t\in[0,T], x\in\mR^d, \mu\in\cM_2(\mR^d), u\in\mU_0.
\end{array}
\label{eq2}
\right.
\ee
 \et 
\begin{proof}
First, we prove sufficiency. For $V\in C_b^{1,2,2}([0,T]\times\mR^d\times\cM_2(\mR^d))$, based on Proposition \ref{itoii}, it holds that
\be
\dif V(t,X_t,\cL_{X_t})&=&(\partial_t+{\bf L}_{b,\sigma,f})V(t,X_t,\cL_{X_t})\dif t+\<(\sigma^*\partial_x V)(t,X_t,\cL_{X_t}), \dif B_t\>\no\\
&&+\int_{\mU_0}\left[V(t, X_{t-}+f(t,X_{t-},\cL_{X_t},u),\cL_{X_t})-V(t, X_{t-},\cL_{X_t})\right]\tilde{N}(\dif t, \dif u).\no\\
\label{ito2}
\ee
Inserting (\ref{eq2}) in (\ref{ito2}), we have 
\ce
\dif V(t,X_t,\cL_{X_t})&=&g_1(t,X_t,\cL_{X_t})\dif t+\int_{\mU_0}g_4(t,X_t,\cL_{X_t},u)\nu(\dif u)\dif t+\<g_2(t,X_t,\cL_{X_t}),\dif B_t\>\\
&&+\int_{\mU_0}g_3(t,X_{t-},\cL_{X_t},u)\tilde{N}(\dif t, \dif u).
\de
By integrating the above equality from $s$ to $t$, one can obtain (\ref{defi}). That is, $F_{s,t}$ is path independent.

Next, let us show necessity. On one hand, since $F_{s,t}$ is path independent, it follows from Definition \ref{pathinde} that 
\be
V(t,X_t,\cL_{X_t})-V(0,X_0,\cL_{X_0})&=&\int_0^t g_1(r,X_r, \cL_{X_r})\dif r+\int_0^t \<g_2(r,X_r, \cL_{X_r}), \dif B_r\>\no\\
&&+\int_0^t\int_{\mU_0}g_3(r,X_{r-}, \cL_{X_r}, u)\tilde{N}(\dif r, \dif u)\no\\
&&+\int_0^t\int_{\mU_0}g_4(r,X_r, \cL_{X_r}, u)\nu(\dif u)\dif r, \quad t\geq 0.
\label{pathinde2}
\ee
On the other hand, by integrating (\ref{ito2}) from $0$ to $t$, we get that 
\be
V(t,X_t,\cL_{X_t})-V(0,X_0,\cL_{X_0})&=&\int_0^t(\partial_r+{\bf L}_{b,\sigma,f})V(r,X_r,\cL_{X_r})\dif r\no\\
&&+\int_0^t\<(\sigma^*\partial_x V)(r,X_r,\cL_{X_r}), \dif B_r\>\no\\
&&+\int_0^t\int_{\mU_0}\bigg[V(r, X_{r-}+f(r,X_{r-},\cL_{X_r},u),\cL_{X_r})\no\\
&&\qquad\qquad -V(r, X_{r-},\cL_{X_r})\bigg]\tilde{N}(\dif r, \dif u).
\label{ito3}
\ee
Thus, $V(t,X_t,\cL_{X_t})-V(0,X_0,\cL_{X_0})$ has two expressions. Since $V(t,X_t,\cL_{X_t})-V(0,X_0,\cL_{X_0})$ is a semimartingale, by uniqueness for decomposition of the semimartingale it holds that
\ce
&&g_1(r,X_r, \cL_{X_r})+\int_{\mU_0}g_4(r,X_r, \cL_{X_r}, u)\nu(\dif u)=(\partial_r+{\bf L}_{b,\sigma,f})V(r,X_r,\cL_{X_r}),\\
&&g_2(r,X_r, \cL_{X_r})=(\sigma^*\partial_x V)(r,X_r,\cL_{X_r}),\\
&&g_3(r,X_r, \cL_{X_r}, u)=V(r, X_r+f(r,X_r,\cL_{X_r},u),\cL_{X_r})-V(r, X_r,\cL_{X_r}), r\in[0,T].
\de
And then for any $s\in[0,T]$ and $\mu=\cL_{X_s}\in\cM_2(\mR^d)$, we know that
\ce
&&g_1(s,X_s, \mu)+\int_{\mU_0}g_4(s,X_s, \mu, u)\nu(\dif u)=(\partial_r+{\bf L}_{b,\sigma,f})V(s,X_s,\mu),\\
&&g_2(s,X_s, \mu)=(\sigma^*\partial_x V)(s,X_s,\mu),\\
&&g_3(s,X_s, \mu, u)=V(s, X_s+f(s,X_s,\mu,u),\mu)-V(s, X_s,\mu),
\de
and then
\ce
&&g_1(s,x, \mu)+\int_{\mU_0}g_4(s,x, \mu, u)\nu(\dif u)=(\partial_r+{\bf L}_{b,\sigma,f})V(s,x,\mu),\\
&&g_2(s,x, \mu)=(\sigma^*\partial_x V)(s,x,\mu),\\
&&g_3(s,x, \mu, u)=V(s, x+f(s,x,\mu,u),\mu)-V(s, x,\mu), x\in {\rm supp}(\mu).
\de

To show (\ref{eq2}), we replace $\mu$ by $\mu^n=\mu\ast\cN_d(0, \frac{1}{n}I_d)$ in the above equality, where $\cN_d(0, \frac{1}{n}I_d)$ denotes the $d$-dimensional Gaussian distribution with the mean $0$ and the covariance matrix $\frac{1}{n}I_d$. Note that ${\rm supp}(\cN_d(0, \frac{1}{n}I_d))=\mR^d$, ${\rm supp}(\mu^n)=\mR^d$ and $\mu^n\rightarrow\mu$ as $n\rightarrow\infty$, which together with continuity of all the related functions in $\mu$ yields (\ref{eq2}). The proof is complete.
\end{proof}

\subsection{Some special cases}
In the subsection, we analyze some special cases.

First of all, we give out a solution of the partial integro-differential equation in (\ref{eq2}). To do this, we introduce two McKean-Vlasov stochastic differential equations with jumps on $\mR^d$: for any $\xi\in L^2(\Omega, \sF, \mP; \mR^d)$ and $x\in\mR^d$,
\be
X_t^{s,\xi}&=&\xi+\int_s^t b(r,X^{s,\xi}_r,\cL_{X^{s,\xi}_r})\dif r+\int_s^t\sigma(r,X^{s,\xi}_r,\cL_{X^{s,\xi}_r})\dif B_r\no\\
&&+\int_s^t\int_{\mU_0}f(r,X^{s,\xi}_{r-},\cL_{X^{s,\xi}_r},u)\tilde{N}(\dif r, \dif u), \quad 0\leq s<t\leq T,\label{mv1}\\
X_t^{s,x,\xi}&=&x+\int_s^t b(r,X^{s,x,\xi}_r,\cL_{X^{s,\xi}_r})\dif r+\int_s^t\sigma(r,X^{s,x,\xi}_r,\cL_{X^{s,\xi}_r})\dif B_r\no\\
&&+\int_s^t\int_{\mU_0}f(r,X^{s,x,\xi}_{r-},\cL_{X^{s,\xi}_r},u)\tilde{N}(\dif r, \dif u), \quad 0\leq s<t\leq T,
\label{mv2}
\ee
and a backward McKean-Vlasov stochastic differential equation
\be\left\{\begin{array}{ll}
\dif Y_t^{s,x,\xi}=g_1(t, X_t^{s,x,\xi},\cL_{X_t^{s,\xi}})\dif t+\int_{\mU_0}g_4(t, X_t^{s,x,\xi},\cL_{X_t^{s,\xi}},u)\nu(\dif u)\dif t,\\
Y_T^{s,x,\xi}=\Phi(X_T^{s,x,\xi},\cL_{X_T^{s,\xi}}).
\end{array}
\right.
\label{mv3}
\ee
If $b,\sigma$ are bounded, $|f(t,x,\mu,u)|\leq C\|u\|_{\mU}$ for $(t,x,\mu,u)\in[0,T]\times\mR^d\times\cM_2(\mR^d)\times\mU_0$ and some constant $C>0$,  and $K_1(\cdot), K_2(\cdot)$ are two constants, under (${\bf H}_{b,\sigma}$) (${\bf H}_f$), based on \cite[Theorem 3.1, Page 7]{HL}, it holds that the above equations (\ref{mv1}) (\ref{mv2}) have unique solutions $X_t^{s,\xi}, X_t^{s,x,\xi}$, respectively. If we further assume $g_1, \int_{\mU_0}g_4(\cdot, \cdot, \cdot, u)\nu(\dif u), \Phi$ are bounded, the above equation (\ref{mv3}) also has a unique solution. For $\Phi\in C_b^{2,2;1}(\mR^d\times\cM_2(\mR^d)), g_1(t, \cdot, \cdot)\in C_b^{2,2;1}(\mR^d\times\cM_2(\mR^d)), \int_{\mU_0}g_4(t, \cdot, \cdot, u)\nu(\dif u)\in C_b^{2,2;1}(\mR^d\times\cM_2(\mR^d))$ and $g_1(\cdot, x, \mu)\in C([0,T]),  \int_{\mU_0}g_4(\cdot, x, \mu, u)\nu(\dif u)\in C([0,T])$, set
\be
V(t,x,\mu)&:=&\mE\bigg[\Phi(X_T^{t,x,\xi},\cL_{X_T^{t,\xi}})-\int_t^Tg_1(r, X_r^{t,x,\xi},\cL_{X_r^{t,\xi}})\dif r\no\\
&&\quad -\int_t^T\int_{\mU_0}g_4(r, X_r^{t,x,\xi},\cL_{X_r^{t,\xi}},u)\nu(\dif u)\dif r\bigg], \quad \mu=\cL_{\xi},
\label{vdef1}
\ee
and  then by \cite[Theorem 9.2, Page 3159]{L}, it holds that $V(t,x,\mu)\in C_b^{1,2,2}([0,T]\times\mR^d\times\cM_2(\mR^d))$ is the unique solution of the following nonlocal partial integral-differential equation
\ce\left\{\begin{array}{ll}
(\partial_t+{\bf L}_{b,\sigma,f})V(t,x,\mu)=g_1(t, x,\mu)+\int_{\mU_0}g_4(t, x,\mu,u)\nu(\dif u), \\
V(T, x,\mu)=\Phi(x,\mu),\\
t\in[0,T], x\in\mR^d, \mu\in\cM_2(\mR^d).
\end{array}
\right.
\de
Thus, by combining Theorem \ref{funthe} with \cite[Theorem 9.2, Page 3159]{L}, one can have the following result.

\bc\label{concfunc}
Assume that (${\bf H}_{b,\sigma}$) (${\bf H}_f$) hold, $(b(t, x, \mu), \sigma(t, x, \mu))\in C_b^{1,2,2;1}([0,T]\times\mR^d\times\cM_2(\mR^d)\mapsto\mR^d\times\mR^{d\times m})$, and all the derivatives of $f(t, x, \mu, u)$ in $t$ order $1$ and in $x, \mu$ up to order $2$ are bounded by $L\|u\|_{\mU}$ and Lipschitz continuous with a Lipschitz factor $L\|u\|_{\mU}$. Then for $V(t,x,\mu)$ defined in (\ref{vdef1}), $g_2\in C([0,T]\times\mR^d\times\cM_2(\mR^d)\mapsto\mR^m)$ and $g_3(\cdot,\cdot,\cdot,u)\in C([0,T]\times\mR^d\times\cM_2(\mR^d)\mapsto\mR)$, $F_{s, t}$ is path independent in the sense of (\ref{defi}) if and only if $V, g_2, g_3$ satisfy
\ce\left\{\begin{array}{ll}
(\sigma^*\partial_x V)(t, x,\mu)=g_2(t, x,\mu),\\
V\(t, x+f(t,x,\mu,u),\mu\)-V(t, x,\mu)=g_3(t, x,\mu,u),\\
 t\in[0,T], x\in\mR^d, \mu\in\cM_2(\mR^d), u\in\mU_0.
\end{array}
\right.
\de
\ec

In the following we analyze some special cases of $F_{s,t}$. If $g_1=0, g_4=0$, $F_{s, t}$ reduces to 
$$
F^{g_2, g_3}_{s,t}:=\int_s^t \<g_2(r,X_r, \cL_{X_r}), \dif B_r\>+\int_s^t\int_{\mU_0}g_3(r,X_{r-}, \cL_{X_r}, u)\tilde{N}(\dif r, \dif u).
$$
We follow up the above deduction to get that for $V(t,x,\mu):=\mE[\Phi(X_T^{t,x,\xi},\cL_{X_T^{t,\xi}})]$, $g_2\in C([0,T]\times\mR^d\times\cM_2(\mR^d)\mapsto\mR^m)$ and $g_3(\cdot,\cdot,\cdot,u)\in C([0,T]\times\mR^d\times\cM_2(\mR^d)\mapsto\mR)$, $F^{g_2, g_3}_{s, t}$ is path independent in the sense of (\ref{defi}) if and only if $V, g_2, g_3$ satisfy
\ce\left\{\begin{array}{ll}
(\sigma^*\partial_x V)(t, x,\mu)=g_2(t, x,\mu),\\
V\(t, x+f(t,x,\mu,u),\mu\)-V(t, x,\mu)=g_3(t, x,\mu,u),\\
 t\in[0,T], x\in\mR^d, \mu\in\cM_2(\mR^d), u\in\mU_0.
\end{array}
\right.
\de
This result also can be obtained by Theorem \ref{funthe} and \cite[Theorem 7.3, Page 47]{HL}. If $g_1=\frac{1}{2\beta}|g_2|^2, \beta\neq0$, $F_{s, t}$ reduces to 
\ce
F^{g_2,g_3,g_4}_{s,t}&:=&\int_s^t \frac{1}{2\beta}|g_2|^2(r,X_r, \cL_{X_r})\dif r
+\int_s^t \<g_2(r,X_r, \cL_{X_r}), \dif B_r\>\no\\
&&+\int_s^t\int_{\mU_0}g_3(r,X_{r-}, \cL_{X_r}, u)\tilde{N}(\dif r, \dif u)\\
&&+\int_s^t\int_{\mU_0}g_4(r,X_r, \cL_{X_r}, u)\nu(\dif u)\dif r, \quad 0\leq s< t\leq T.
\de
Thus, by Theorem \ref{funthe}, it holds that $F^{g_2,g_3,g_4}_{s,t}$ is path independent in the sense of (\ref{defi}) if and only if 
\be\left\{\begin{array}{ll}
(\partial_t+{\bf L}_{b,\sigma,f})V(t,x,\mu)=\frac{1}{2\beta}|\sigma^*\partial_x V|^2(t, x,\mu)+\int_{\mU_0}g_4(t, x,\mu,u)\nu(\dif u),\\
(\sigma^*\partial_x V)(t, x,\mu)=g_2(t, x,\mu),\\
V\(t, x+f(t,x,\mu,u),\mu\)-V(t, x,\mu)=g_3(t, x,\mu,u), \\
t\in[0,T], x\in\mR^d, \mu\in\cM_2(\mR^d), u\in\mU_0.
\end{array}
\label{eq5}
\right.
\ee

\subsection{The relationship between Theorem \ref{funthe} and some known results}
In the subsection, we discuss the relationship between Theorem \ref{funthe} and \cite[Theorem 2.6]{qw1} \cite[Theorem 2.2]{rw}. 

Let $\lambda: [0,\infty)\times \mU_0\to(0,1]$ be a measurable function. Here we require that if $\lambda(t,u)=1$, $f(t,x,\mu,u)=0$ immediately. And then there exists an integer-valued $(\sF_t)_{t\geq 0}$-Poisson
random measure $N_{\lambda}(\dif t, \dif u)$ on $(\Omega,\sF,\mP;(\sF_t)_{t\geq 0})$ with intensity
$\mE(N_{\lambda}(\dif t, \dif u))=\lambda(t,u)\dif t\nu(\dif u)$. Denote
$$\tilde{N}_\lambda(\dif t,\dif u):=N_\lambda(\dif t,\dif u)-\lambda(t,u)\dif t\nu(\dif u)$$
that is, $\tilde{N}_\lambda(\dif t,\dif u)$ stands for the compensated $(\sF_t)_{t\geq 0}$-predictable martingale measure of $N_{\lambda}(\dif t, \dif u)$. Moreover, $\tilde{N}_\lambda(\dif t,\dif u)$ is independent of $B_{\cdot}$. We replace $\tilde{N}(\dif t, \dif u)$ by $\tilde{N}_\lambda(\dif t,\dif u)$ in Eq.(\ref{Eq1}). Thus, the solution of the new equation is denoted as $X_t^{\lambda}$. Besides, we assume that there exists a measurable function $\tilde{b}: [0,T]\times\mR^d\times\cM_2(\mR^d)\rightarrow\mR^m$ such that $b=\sigma \tilde{b}$. For convenience of the following deduction, we also assume:
\ce
&&\mE\Big[\exp\Big\{\frac{1}{2}\int_0^T\left|\tilde{b}(s,X^{\lambda}_s, \cL_{X^{\lambda}_s})\right|^2\dif s\Big\}\Big]<\infty,\\
&& \int_0^T\int_{\mU_0}\left(\frac{1-\lambda(s,u)}{\lambda(s,u)}\right)^2\lambda(s,u)\nu(\dif u)\dif s
<\infty.
\de
So, set 
\ce
\Gamma_t:&=&\exp\bigg\{-\int_0^t\<\tilde{b}(s,X^{\lambda}_s, \cL_{X^{\lambda}_s}),\dif B_s\>-\frac{1}{2}\int_0^t
\left|\tilde{b}(s,X^{\lambda}_s, \cL_{X^{\lambda}_s})\right|^2\dif s\\
&&\qquad\quad -\int_0^t\int_{\mU_0}\log\lambda(s,u)\tilde{N}_{\lambda}(\dif s, \dif u)\\
&&\qquad\quad -\int_0^t\int_{\mU_0}\Big(\big(\log\lambda(s,u)\big)\lambda(s,u)+\big(1-\lambda(s,u)\big)\Big)\nu(\dif u)\dif s\bigg\},
\de
and then by the same deduction to that in \cite{qw1}, it holds that $\Gamma_t$ is an exponential martingale. Define a new probability $\mQ$ as
$$
\frac{\dif\mQ}{\dif\mP}=\Gamma_T.
$$
Thus, under $\mQ$, 
$$
\tilde{B}_t:=B_t+\int_0^t\tilde{b}(s,X^{\lambda}_s, \cL_{X^{\lambda}_s})\dif s
$$
is a $d$-dimensional Brownian motion and 
$$
\bar{N}(\dif t,\dif u)=N_\lambda(\dif t,\dif u)-\dif t\nu(\dif u)
$$
is the compensated $(\sF_t)_{t\geq 0}$-predictable martingale measure of $N_{\lambda}(\dif t, \dif u)$. Moreover, Eq.(\ref{Eq1}) becomes 
\ce
X^{\lambda}_t=X^{\lambda}_s+\int_s^t\sigma(r,X^{\lambda}_r,\cL_{X^{\lambda}_r})\dif \tilde{B}_r+\int_s^t\int_{\mU_0}f(r,X^{\lambda}_{r-},\cL_{X^{\lambda}_r},u)\bar{N}(\dif r, \dif u).
\de
That is, $X^{\lambda}_t$ is a local martingale.

Now, take 
\ce
&&g_1(t,x,\mu)=\frac{1}{2}
\left|\tilde{b}(t,x,\mu)\right|^2, ~ g_2(t,x,\mu)=\tilde{b}(t,x,\mu),\\
&&g_3(t,x,\mu,u)=\log\lambda(t,u), \quad g_4(t,x,\mu,u)=\log\lambda(t,u)+\left(\frac{1}{\lambda(t,u)}-1\right).
\de
And then by Theorem \ref{funthe}, we know that $F^{\lambda}_{0,t}:=-\log\Gamma_t$ is path independent in the sense of (\ref{defi}) if and only if 
\be\left\{\begin{array}{ll}
(\partial_t+{\bf L}_{b,\sigma,f})V(t,x,\mu)\\
=\frac{1}{2}\left|(\sigma^*\partial_x V)(t,x,\mu)\right|^2+\int_{\mU_0}\Big(\big(\log\lambda(t,u)\big)\lambda(t,u)+\big(1-\lambda(t,u)\big)\Big)\nu(\dif u),\\
(\sigma^*\partial_x V)(t, x,\mu)=\tilde{b}(t,x,\mu),\\
V\(t, x+f(t,x,\mu,u),\mu\)-V(t, x,\mu)=\log\lambda(t,u), \\
t\in[0,T], x\in\mR^d, \mu\in\cM_2(\mR^d), u\in\mU_0.
\end{array}
\right.
\label{eq6}
\ee
The equation is just right Eq.(\ref{eq5}) with $\beta=1$. And then we rewrite Eq.(\ref{eq6}) as 
\ce\left\{\begin{array}{ll}
\partial_tV(t,x,\mu)={\bf L}_{\sigma,f}V(t,x,\mu),\\
b(t,x,\mu)=(\sigma\sigma^*\partial_x V)(t, x,\mu),\\
\lambda(t,u)=\exp\left\{V\(t, x+f(t,x,\mu,u),\mu\)-V(t, x,\mu)\right\}, \\
t\in[0,T], x\in\mR^d, \mu\in\cM_2(\mR^d), u\in\mU_0,
\end{array}
\right.
\de
where 
\ce
{\bf L}_{\sigma,f}V(t,x,\mu)&:=&-\frac{1}{2}\tr\((\sigma\sigma^*)\partial_x^2 V\)(t, x,\mu)-\frac{1}{2}\left|(\sigma^*\partial_x V)(t,x,\mu)\right|^2\no\\
&&-\int_{\mR^d}\<(\sigma\sigma^*\partial_y V)(t, y,\mu), (\partial_\mu V)(t, x,\mu)(y)\>\mu(\dif y)\no\\
&&-\frac{1}{2}\int_{\mR^d}\tr\((\sigma\sigma^*)(t, y,\mu)\partial_y\partial_\mu V(t, x,\mu)(y)\)\mu(\dif y) \no\\
&&-\int_{\mU_0}\bigg[e^{V(t, x+f(t,x,\mu,u),\mu)-V(t, x,\mu)}-1\no\\
&&\qquad\qquad\qquad -\<f(t,x,\mu,u), \partial_x V(t, x,\mu)\>e^{V(t, x+f(t,x,\mu,u),\mu)-V(t, x,\mu)}\bigg]\nu(\dif u)\no\\
&&-\int_{\mU_0}\int_0^1\int_{\mR^d}\bigg<\partial_\mu V(t,x,\mu)\(y+\eta f(t,y,\mu,u)\)\no\\
&&\qquad\qquad -\partial_\mu V(t,x,\mu)(y), f(t,y,\mu,u)\bigg>\mu(\dif y)\dif \eta \nu(\dif u).
\de
If $b, \sigma, f, V$ are independent of $\mu$, this is just right Theorem 2.6 in \cite{qw1}. Moreover, if $f=0$,  this is exactly Theorem 2.2 in \cite{rw}. Therefore, our result is more general.

\section{Appendix}\label{app}
{\bf The proof of Proposition \ref{itoii}.}

Set $\mu_t:=\cL_{X_t}$, 
and then $h(t,X_t,\cL_{X_t})=h(t,X_t,\mu_t)$. Define $\bar{h}(t,x):=h(t,x,\mu_t)$, and then $\bar{h}(t,X_t)=h(t,X_t,\mu_t)$. Moreover, based on $h\in C_b^{1,2,2}([0,T]\times\mR^d\times\cM_2(\mR^d))$, it holds that $\bar{h}$ is $C^2$ in $x$. However, we don't know the differentiability of $\bar{h}$ in $t$. Note that the differentiability of $\bar{h}$ in $t$ comes from two parts-$h(t,x,\mu)$ in $t$ for fixed $x,\mu$ and $h(s,x,\mu_t)$ in $t$ for fixed $s,x$. Therefore, to apply the classical It\^o formula to $\bar{h}(t,X_t)$, we only need to consider the second part by $h\in C_b^{1,2,2}([0,T]\times\mR^d\times\cM_2(\mR^d))$.

{\bf Step 1.} Assume that $b,\sigma$ are bounded, $|f(t,x,\mu,u)|\leq C\|u\|_{\mU}$ for $(t,x,\mu,u)\in[0,T]\times\mR^d\times\cM_2(\mR^d)\times\mU_0$ and some constant $C>0$. We study the differentiability of $h(s,x,\mu_t)$ in $t$.

Here, we follow the method in \cite{ccd} to deal with it. For the convenience to our expression, we take $H(\mu_t):=h(s,x,\mu_t)$. For any positive integer $K$, set
\be
x^1, x^2,\cdots, x^K\in\mR^d, H^K(x^1, x^2,\cdots, x^K):=H\(\frac{1}{K}\sum_{l=1}^K\delta_{x^l}\),
\label{hkde}
\ee
and then $H^K(x^1, x^2,\cdots, x^K)$ is a function on $\mR^{d\times K}$. Moreover, by \cite[Proposition 3.1, Page 15]{ccd}, it holds that $H^K$ is $C^2$ on $\mR^{d\times K}$ and 
\be
\partial_{x^i}H^K(x^1, x^2,\cdots, x^K)&=&\frac{1}{K}\partial_{\mu}H\(\frac{1}{K}\sum_{l=1}^K\delta_{x^l}\)(x^i),\no\\
\partial^2_{x^ix^j}H^K(x^1, x^2,\cdots, x^K)&=&\frac{1}{K}\partial_{y}\partial_{\mu}H\(\frac{1}{K}\sum_{l=1}^K\delta_{x^l}\)(x^i)\delta_{i,j}+\frac{1}{K^2}\partial^2_{\mu}H\(\frac{1}{K}\sum_{l=1}^K\delta_{x^l}\)(x^i, x^j),\no\\
\qquad\qquad i,j=1,2,\cdots, K,
\label{deri}
\ee
where $\delta_{i,j}=1, i=j$, $\delta_{i,j}=0, i\neq j$.
Besides, we take $K$ independent copies $X^l_t, l=1,2,\cdots, K$ of $X_t$. That is, 
$$
\dif X^l_t=b(t,X^l_t,\cL_{X^l_t})\dif t+\sigma(t,X^l_t,\cL_{X^l_t})\dif B^l_t+\int_{\mU_0}f(t,X^l_{t-},\cL_{X^l_t},u)\tilde{N}^l(\dif t, \dif u), \quad l=1,2,\cdots, K,
$$
where $B^l, N^l, l=1,2,\cdots, K$ are mutually independent and have the same distributions to that of $B, N$, respectively. And then applying the It\^o formula to $H^K(X_t^1, X_t^2, \cdots, X_t^K)$ and taking the expectation on both sides, we obtain that for $0\leq t<t+v\leq T$
\ce
&&\mE H^K(X_{t+v}^1, X_{t+v}^2, \cdots, X_{t+v}^K)\\
&=&\mE H^K(X_t^1, X_t^2, \cdots, X_t^K)+\sum_{i=1}^K\int_t^{t+v}\mE\partial_{x^i}H^K(X_s^1, X_s^2, \cdots, X_s^K)b(s,X^i_s,\cL_{X^i_s})\dif s\\
&&+\frac{1}{2}\sum_{i=1}^K\int_t^{t+v}\mE\partial^2_{x^ix^i}H^K(X_s^1, X_s^2, \cdots, X_s^K)\sigma\sigma^*(s,X^i_s,\cL_{X^i_s})\dif s\\
&&+\int_t^{t+v}\int_{\mU_0}\mE\[H^K\(X_s^1+f(s,X^1_{s},\cL_{X^1_s},u), X_s^2, \cdots, X_s^K\)-H^K(X_s^1, X_s^2, \cdots, X_s^K)\\
&&\qquad-\partial_{x^1}H^K(X_s^1, X_s^2, \cdots, X_s^K)f(s,X^1_{s},\cL_{X^1_s},u)\]\nu(\dif u)\dif s\\
&&+\cdots\\
&&+\int_t^{t+v}\int_{\mU_0}\mE\[H^K\(X_s^1, X_s^2, \cdots, X_s^K+f(s,X^K_{s},\cL_{X^K_s},u)\)-H^K(X_s^1, X_s^2, \cdots, X_s^K)\\
&&\qquad-\partial_{x^K}H^K(X_s^1, X_s^2, \cdots, X_s^K)f(s,X^K_{s},\cL_{X^K_s},u)\]\nu(\dif u)\dif s\\
&=&\mE H^K(X_t^1, X_t^2, \cdots, X_t^K)+K\int_t^{t+v}\mE\partial_{x^1}H^K(X_s^1, X_s^2, \cdots, X_s^K)b(s,X^1_s,\cL_{X^1_s})\dif s\\
&&+\frac{K}{2}\int_t^{t+v}\mE\partial^2_{x^1x^1}H^K(X_s^1, X_s^2, \cdots, X_s^K)\sigma\sigma^*(s,X^1_s,\cL_{X^1_s})\dif s\\
&&+K\int_t^{t+v}\int_{\mU_0}\int_0^1\mE\[\(\partial_{x^1}H^K(X_s^1+\eta f(s,X^1_{s},\cL_{X^1_s},u), X_s^2, \cdots, X_s^K)\\
&&\qquad\qquad -\partial_{x^1}H^K(X_s^1, X_s^2, \cdots, X_s^K)\)f(s,X^1_{s},\cL_{X^1_s},u)\]\dif\eta\nu(\dif u)\dif s,
\de
where the property of the same distributions for $X^l_t, l=1,2,\cdots, K$ is used in the second equality. Inserting (\ref{hkde}) (\ref{deri}) in the above equality, we get that
\ce
&&\mE H\(\frac{1}{K}\sum_{l=1}^K\delta_{X_{t+v}^l}\)\\
&=&\mE H\(\frac{1}{K}\sum_{l=1}^K\delta_{X_t^l}\)+\int_t^{t+v}\mE\partial_{\mu}H\(\frac{1}{K}\sum_{l=1}^K\delta_{X_s^l}\)(X_s^1)b(s,X^1_s,\cL_{X^1_s})\dif s\\
&&+\frac{1}{2}\int_t^{t+v}\mE\partial_y\partial_{\mu}H\(\frac{1}{K}\sum_{l=1}^K\delta_{X_s^l}\)(X_s^1)\sigma\sigma^*(s,X^1_s,\cL_{X^1_s})\dif s\\
&&+\frac{1}{2K}\int_t^{t+v}\mE\partial^2_{\mu}H\(\frac{1}{K}\sum_{l=1}^K\delta_{X_s^l}\)(X_s^1,X_s^1)\sigma\sigma^*(s,X^1_s,\cL_{X^1_s})\dif s\\
&&+\int_t^{t+v}\int_{\mU_0}\int_0^1\mE\[\(\partial_{\mu}H\(\frac{1}{K}\delta_{X_s^1+\eta f(s,X^1_{s},\cL_{X^1_s},u)}+\frac{1}{K}\sum_{l=2}^K\delta_{X_s^l}\)\circ(X_s^1+\eta f(s,X^1_{s},\cL_{X^1_s},u))\\
&&\qquad\qquad\qquad-\partial_{\mu}H\(\frac{1}{K}\sum_{l=1}^K\delta_{X_s^l}\)(X_s^1)\)f(s,X^1_{s},\cL_{X^1_s},u)\]\dif\eta\nu(\dif u)\dif s.
\de

Next, we take the limit on both sides of the above equality. Note that 
\ce
\lim_{K\rightarrow\infty}\mE\left[\sup_{0\leq t\leq T}\rho^2\left(\frac{1}{K}\sum_{l=1}^K\delta_{X_t^l},\mu_t\right)\right]=0,
\de
which comes from \cite[Section 5]{HK}. And then as $K\rightarrow\infty$, by continuity and boundedness of $H, \partial_{\mu}H, \partial_y\partial_{\mu}H$, and boundedness of $\partial^2_{\mu}H, b, \sigma$, it follows from the dominated convergence theorem that
\ce
H(\mu_{t+v})&=&H(\mu_t)+\int_t^{t+v}\mE\partial_{\mu}H\(\mu_s\)(X_s^1)b(s,X^1_s,\cL_{X^1_s})\dif s\\
&&+\frac{1}{2}\int_t^{t+v}\mE\partial_y\partial_{\mu}H\(\mu_s\)(X_s^1)\sigma\sigma^*(s,X^1_s,\cL_{X^1_s})\dif s\\
&&+\int_t^{t+v}\int_{\mU_0}\int_0^1\mE\[\(\partial_{\mu}H(\mu_s)\(X_s^1+\eta f(s,X^1_{s},\cL_{X^1_s},u)\)\\
&&\qquad\qquad\qquad-\partial_{\mu}H(\mu_s)(X_s^1)\)f(s,X^1_{s},\cL_{X^1_s},u)\]\dif\eta\nu(\dif u)\dif s.
\de
Thus, by simple calculus we obtain that 
\be
&&\partial_tH(\mu_t)\no\\
&=&\int_{\mR^d}\<b(t, y,\mu_t), \partial_\mu H(\mu_t)(y)\>\mu_t(\dif y)+\frac{1}{2}\int_{\mR^d}\tr\((\sigma\sigma^*)(t, y,\mu_t)\partial_y\partial_\mu H(\mu_t)(y)\)\mu_t(\dif y)\no\\
&&+\int_{\mU_0}\int_0^1\int_{\mR^d}\[\(\partial_{\mu}H(\mu_t)\(y+\eta f(t,y,\mu_t,u)\)-\partial_{\mu}H(\mu_t)(y)\)f(t,y,\mu_t,u)\]\mu_t(\dif y)\dif\eta\nu(\dif u).\no\\
\label{diff}
\ee

{\bf Step 2.} Assume that (${\bf H}_{b,\sigma}$) (${\bf H}_f$) hold. We prove the differentiability of $h(s,x,\mu_t)$ in $t$.

First of all, we choose a smooth function $\chi_n: \mR^d\mapsto\mR^d$ satisfying $\chi_n(x)=x, |x|\leq n$ and $\chi_n(x)=0, |x|\geq 2n$ such that for $x\in\mR^d$
\be
|\chi_n(x)|\leq C, \quad \|\partial \chi_n(x)\|\leq C,
\label{chpro}
\ee
where the positive constant $C$ is independent of $n$. Set 
$$
b^{(n)}(t,x,\mu):=b(t,\chi_n(x),\mu), \sigma^{(n)}(t,x,\mu):=\sigma(t,\chi_n(x),\mu), f^{(n)}(t,x,\mu,u):=f(t,\chi_n(x),\mu,u),
$$
and then 
$$
b^{(n)}(t,x,\mu)\rightarrow b(t,x,\mu), \sigma^{(n)}(t,x,\mu)\rightarrow \sigma(t,x,\mu), f^{(n)}(t,x,\mu,u)\rightarrow f(t,x,\mu,u)
$$ 
as $n\rightarrow \infty$. Moreover, by Remark \ref{lingro} we know that $b^{(n)},\sigma^{(n)}$ are bounded, $|f^{(n)}(t,x,\mu,u)|\leq C\|u\|_{\mU}$ for $(t,x,\mu,u)\in[0,T]\times\mR^d\times\cM_2(\mR^d)\times\mU_0$, and $b^{(n)},\sigma^{(n)}, f^{(n)}$ satisfy (${\bf H}_{b,\sigma}$) (${\bf H}_f$). Thus, the following equation
\ce 
\dif X^{(n)}_t=b^{(n)}(t,X^{(n)}_t,\cL_{X^{(n)}_t})\dif t+\sigma^{(n)}(t,X^{(n)}_t,\cL_{X^{(n)}_t})\dif B_t+\int_{\mU_0}f^{(n)}(t,X^{(n)}_{t-},\cL_{X^{(n)}_t},u)\tilde{N}(\dif t, \dif u)
\de
has a unique solution denoted by $X^{(n)}_{\cdot}$ and $\mu_t^{(n)}:=\cL_{X^{(n)}_t}$. And then by {\bf Step 1.}, it holds that for $0\leq t<t+v\leq T$
\be
H(\mu^{(n)}_{t+v})-H(\mu^{(n)}_t)&=&\int_t^{t+v}\int_{\mR^d}\<b^{(n)}(r,y,\mu^{(n)}_r), \partial_\mu H(\mu^{(n)}_r)(y)\>\mu^{(n)}_r(\dif y)\dif r\no\\
&&+\frac{1}{2}\int_t^{t+v}\int_{\mR^d}\tr\((\sigma^{(n)}\sigma^{(n)*})(r, y,\mu^{(n)}_r)\partial_y\partial_\mu H(\mu^{(n)}_r)(y)\)\mu^{(n)}_r(\dif y)\dif r\no\\
&&+\int_t^{t+v}\int_{\mU_0}\int_0^1\int_{\mR^d}\[\(\partial_{\mu}H(\mu^{(n)}_r)\(y+\eta f^{(n)}(r,y,\mu^{(n)}_r,u)\)\no\\
&&\qquad\quad -\partial_{\mu}H(\mu^{(n)}_r)(y)\)f^{(n)}(r,y,\mu^{(n)}_r,u)\]\mu^{(n)}_r(\dif y)\dif\eta\nu(\dif u)\dif r.
\label{bsigcon}
\ee

Next, we observe the limit of $\mu^{(n)}_t$ as $n\rightarrow \infty$ for any $t\in[0,T]$. Applying the It\^o formula to $|X_t^{(n)}-X_t|^2$ and taking the expectation on two sides, one can obtain that
\ce
&&\mE|X_t^{(n)}-X_t|^2\\
&=&2\mE\int_0^t\<b^{(n)}(r,X^{(n)}_r,\mu^{(n)}_r)-b(r,X_r,\mu_r), X_r^{(n)}-X_r\>\dif r\\
&&+\mE\int_0^t\tr\((\sigma^{(n)}(r,X^{(n)}_r,\mu^{(n)}_r)-\sigma(r,X_r,\mu_r))(\sigma^{(n)*}(r,X^{(n)}_r,\mu^{(n)}_r)-\sigma^*(r,X_r,\mu_r))\)\dif r\\
&&+\mE\int_0^t\int_{\mU_0}|f^{(n)}(r,X^{(n)}_r,\mu^{(n)}_r,u)-f(r,X_r,\mu_r,u)|^2\nu(\dif u)\dif r\\
&\leq&C\mE\int_0^t|X_r^{(n)}-X_r|^2\dif r+C\mE\int_0^t\rho^2(\mu^{(n)}_r,\mu_r)\dif r+C\mE\int_0^t|\chi_n(X_r^{(n)})-X_r|^2\dif r\\
&\leq&C\int_0^t\mE|X_r^{(n)}-X_r|^2\dif r+C\mE\int_0^t|\chi_n(X_r)-X_r|^2\dif r,
\de
where we use (${\bf H}_{b,\sigma}$) (${\bf H}_f$) (\ref{chpro}) and $\rho^2(\mu^{(n)}_r,\mu_r)\leq \mE|X_r^{(n)}-X_r|^2$ in the first and second inequalities, respectively. Thus, the Gronwall inequality admits us to have that
\ce
\mE|X_t^{(n)}-X_t|^2\leq C\mE\int_0^T|\chi_n(X_r)-X_r|^2\dif r.
\de
The fact that $|\chi_n(x)|\leq |x|$ for $x\in\mR^d$, together with (\ref{xbou}) and the dominated convergence theorem, yields that
$$
\lim\limits_{n\rightarrow\infty}\mE\int_0^T|\chi_n(X_r)-X_r|^2\dif r=0.
$$
So, we get that
$$
\lim\limits_{n\rightarrow\infty}\sup\limits_{t\in[0,T]}\rho^2(\mu^{(n)}_t,\mu_t)\leq \lim\limits_{n\rightarrow\infty}\sup\limits_{t\in[0,T]}\mE|X_t^{(n)}-X_t|^2=0.
$$
Taking the limit on both sides of (\ref{bsigcon}), by Remark \ref{lingro} and the dominated convergence theorem, one can still obtain (\ref{diff}).

{\bf Step 3.} We prove (\ref{itoeq}).

By {\bf Step 2.}, we know that $\bar{h}(t, x)$ is $C^1$ in $t$ and $C^2$ in $x$. Therefore the classical It\^o formula admits us to obtain that
\ce
&&\dif h(t,X_t,\cL_{X_t})=\dif \bar{h}(t,X_t)\\
&=&\partial_t\bar{h}(t,X_t)\dif t+\<\partial_x \bar{h}(t, X_t), b(t, X_t,\mu_t)\>\dif t+\<\partial_x \bar{h}(t, X_t), \sigma(t,X_t,\mu_t)\dif B_t\>\\
&&+\int_{\mU_0}\bigg[\bar{h}\(t, X_t+f(t,X_t,\mu_t,u)\)-\bar{h}(t, X_t)-\<f(t,X_t,\mu_t,u), \partial_x \bar{h}(t, X_t)\>\bigg]\nu(\dif u)\dif t\\
&&+\int_{\mU_0}\left[\bar{h}\(t, X_{t-}+f(t,X_{t-},\mu_t, u)\)-\bar{h}(t, X_{t-})\right]\tilde{N}(\dif t, \dif u)\\
&&+\frac{1}{2}\tr\(\sigma\sigma^*(t,X_t,\mu_t)\partial^2_x \bar{h}(t, X_t)\)\dif t\\
&=&\partial_t h(t,X_t,\mu_t)\dif t+\partial_t h(s,x,\mu_t)|_{s=t, x=X_t}\dif t+\<\partial_x h(t, X_t, \mu_t), b(t, X_t,\mu_t)\>\dif t\\
&&+\<\partial_x h(t, X_t, \mu_t), \sigma(t,X_t,\mu_t)\dif B_t\>+\frac{1}{2}\tr\(\sigma\sigma^*(t,X_t,\mu_t)\partial^2_x h(t, X_t, \mu_t)\)\dif t\\
&&+\int_{\mU_0}\left[h\(t, X_{t-}+f(t,X_{t-},\mu_t, u),\mu_t\)-h(t, X_{t-},\mu_t)\right]\tilde{N}(\dif t, \dif u)\\
&&+\int_{\mU_0}\bigg[h\(t, X_t+f(t,X_t,\mu_t,u),\mu_t\)-h(t, X_t,\mu_t)\\
&&\qquad\qquad\qquad -\<f(t,X_t,\mu_t,u), \partial_x h(t, X_t,\mu_t)\>\bigg]\nu(\dif u)\dif t\\
&=&\partial_t h(t,X_t,\mu_t)\dif t+\<\partial_x h(t, X_t, \mu_t), b(t, X_t,\mu_t)\>\dif t+\<\partial_x h(t, X_t, \mu_t), \sigma(t,X_t,\mu_t)\dif B_t\>\\
&&+\int_{\mU_0}\left[h\(t, X_{t-}+f(t,X_{t-},\mu_t,u),\mu_t\)-h(t, X_{t-},\mu_t)\right]\tilde{N}(\dif t, \dif u)\\
&&+\frac{1}{2}\tr\(\sigma\sigma^*(t,X_t,\mu_t)\partial^2_x h(t, X_t, \mu_t)\)\dif t+\int_{\mR^d}\<b(t, y,\mu_t), \partial_\mu h(t, X_t,\mu_t)(y)\>\mu_t(\dif y)\dif t\\
&&+\frac{1}{2}\int_{\mR^d}\tr\((\sigma\sigma^*)(t, y,\mu_t)\partial_y\partial_\mu h(t, X_t, \mu_t)(y)\)\mu_t(\dif y)\dif t\\
&&+\int_{\mU_0}\int_0^1\int_{\mR^d}\bigg<\partial_\mu h(t,X_t,\mu_t)\(y+\eta f(t, y,\mu_t, u)\)\\
&&\qquad\qquad\qquad-\partial_\mu h(t,X_t,\mu_t)(y), f(t, y,\mu_t, u)\bigg>\mu_t(\dif y)\dif \eta \nu(\dif u)\dif t \\
&&+\int_{\mU_0}\bigg[h\(t, X_t+f(t,X_t,\mu_t,u),\mu_t\)-h(t, X_t,\mu_t)\\
&&\qquad\qquad\qquad-\<f(t,X_t,\mu_t,u), \partial_x h(t, X_t,\mu_t)\>\bigg]\nu(\dif u)\dif t.
\de
This is just right  (\ref{itoeq}). The proof is complete.

\bigskip

\textbf{Acknowledgements:}

The authors are very grateful to Professor Xicheng Zhang for valuable discussions. The first author also thanks Professor Renming Song for providing her an excellent environment to work in the University of Illinois at Urbana-Champaign.

\end{document}